\documentclass[amstex,12pt, amssymb]{article}

\usepackage{mathtext}
\usepackage[cp1251]{inputenc}
\usepackage[T2A]{fontenc}
\usepackage[dvips]{graphicx}
\usepackage{amsmath}
\usepackage{amssymb}
\usepackage{amsxtra}
\usepackage{latexsym}
\usepackage{ifthen}

\newcounter{lemma}[section]

\newcounter{corollary}[section]

\newcounter{remark}[section]

\newcounter{theorem}[section]

\newcounter{proposition}[section]

\numberwithin{equation}{section}

\def\XXint#1#2#3{{\setbox0=\hbox{$#1{#2#3}{\int}$}
     \vcenter{\hbox{$#2#3$}}\kern-.5\wd0}}

\def\cc{\setcounter{equation}{0}
\setcounter{figure}{0}\setcounter{table}{0}}

\begin{document}

\markboth{\centerline{VLADIMIR RYAZANOV}} {\centerline{BOUNDARY
BEHAVIOR OF CONJUGATE HARMONIC FUNCTIONS}}

\author{{Vladimir Ryazanov}}

\title{{\bf Correlation of boundary behavior\\ of conjugate harmonic functions}}

\maketitle

\large \begin{abstract}

It is established that if a harmonic function $u$ on the unit disk
$\mathbb D$ in $\mathbb C$ has angular limits on a measurable set
$E$ of the unit circle $\partial\mathbb D$, then its conjugate
harmonic function $v$ in $\mathbb D$ also has angular limits a.e. on
$E$ and both boundary functions are measurable on $E$. The result is
extended to arbitrary Jordan domains with rectifiable boundaries in
terms of angular limits and of the natural parameter.
\end{abstract}

%\medskip


\medskip


\large
\cc
\section{Introduction} First of all, recall that a path in $\mathbb D:=\{ z\in\Bbb C : |z|<1\}$
terminating at $\zeta=e^{i\vartheta}\in\partial\mathbb D$ is called
{\bf nontangential} at $\zeta$ if its part in a neighborhood of
$\zeta$ lies inside of an angle in $\mathbb D$ with the vertex at
$\zeta$. Hence limits along all nontangential paths at $\zeta$ are
also named {\bf angular} at $\zeta$. The latter is a traditional
tool of the geometric function theory, see e.g. monographs
\cite{Du}, \cite{GM}, \cite{Ko}, \cite{L}, \cite{Po} and \cite{P}.
Note that every closed rectifiable Jordan curve has a tangent a.e.
with respect to the natural parameter and the angular limit has the
same sense at its points with a tangent.

\medskip

It is known the very delicate fact due to Lusin that harmonic
functions in the unit disk with continuous (even absolutely
continuous !) boundary data can have conjugate harmonic functions
whose boundary data are not con\-ti\-nuous functions, furthermore,
they can fail to be essentially bounded in neighborhoods of each
point of the unit circle, see e.g. Theorem VIII.13.1 in \cite{Bari}.
Thus, a correlation between boundary data of conjugate harmonic
functions is not a simple matter, see also I.E in \cite{Ko}.

\medskip

Denote by $h^p$, $p\in(0,\infty)$, the class of all harmonic
functions $u$ in $\mathbb D$ with
$$
\sup\limits_{r\in(0,1)}\ \left\{\int\limits_{0}\limits^{2\pi}\
|u(re^{i\vartheta})|^p\ d\vartheta\right\}^{\frac{1}{p}}\ <\ \infty\
.
$$
It is clear that $h^p\subseteq h^{p^{\prime}}$ for all
$p>p^{\prime}$ and, in particular, $h^p\subseteq h^{1}$ for all
$p>1$.

\bigskip

{\bf Remark 1.} It is important that every function in the class
$h^1$ has a.e. nontangential boundary limits, see e.g. Corollary
IX.2.2 in \cite{Go}.

\bigskip

It is also known that a harmonic function $u$ in $\Bbb D$ can be
represented as the Poisson integral
\begin{equation}\label{Poisson}
u(re^{i\vartheta})\ =\ \frac{1}{2\pi}\
\int\limits_{0}\limits^{2\pi}\frac{1-r^2}{1-2r\cos(\vartheta-t)+r^2}\
\varphi(t)\ dt
\end{equation} with a function $\varphi\in L^p(0,2\pi),\: p>1$, if and only if $u\in
h^p$, see e.g. Theorem IX.2.3 in \cite{Go}. Thus, $u(z)\to
\varphi(\vartheta)$ as $z\to e^{i\vartheta}$ along any nontangential
path for a.e. $\vartheta$, see e.g. Corollary IX.1.1 in \cite{Go}.
Moreover, $u(z)\to \varphi(\vartheta_0)$ as $z\to e^{i\vartheta_0}$
at points $\vartheta_0$ of continuity of the function $\varphi$, see
e.g. Theorem IX.1.1 in \cite{Go}.

\medskip

Note also that $v\in h^p$ whenever $u\in h^p$ for all $p>1$ by the
M. Riesz theorem, see \cite{R}, see also Theorem IX.2.4 in
\cite{Go}. Generally speaking, this fact is not trivial but it
follows immediately for $p=2$ from the Parseval equality, see e.g.
the proof of Theorem IX.2.4 in \cite{Go}. The case $u\in h^1$ is
more complicated.

\medskip

The correlation of the boundary behavior of conjugate harmonic
functions outside the classes $h^p$ was not investigated at all.
This is just the subject of the present article.

\bigskip

\cc
\section{The case of the unit disk}

Here we apply in a certain part a construction of Luzin-Priwalow
from the proof of their theorem on the boundary uniqueness for
analytic functions, see \cite{LP}, see also \cite{Ko}, Section
III.D.1.

\bigskip

{\bf Theorem 1.} {\it Let $u:\mathbb D\to\mathbb R$ be a harmonic
function that has angular limits on a measurable set $E$ of the unit
circle $\partial\mathbb D$. Then its conjugate harmonic functions
$v$ in $\mathbb D$ also have (finite !) angular limits a.e. on $E$
and both boundary functions are measurable on $E$.}

\bigskip

{\bf Remark 2.} By the Luzin-Priwalow uniqueness theorem for
meromorphic functions $u$ as well as $v$ cannot have infinite
angular limits on a subset of $\partial\mathbb D$ of a positive
measure, see Section IV.2.5 in \cite{P}.

\bigskip

{\it Proof.} By Remark 2 we may consider that angular limits of $u$
are finite eve\-ry\-whe\-re on the set $E$. Moreover, the measurable
set $E$ admits a countable exhaustion by measure of the arc length
with its closed subsets, see e.g. Theorem III(6.6) in \cite{S}, and
hence with no loss of generality we may also consider that $E$ is
compact, see e.g. Proposition I.9.3 in \cite{Bou}.

Following \cite{Ko}, Section III.D.1, we set, for
$\zeta\in\partial\mathbb D$,
\begin{equation}\label{triangle}
S_{\zeta}\ =\ \left\{\ z\in\mathbb D:\ |z|\ >\ \frac{1}{\surd {2}}\
,\ |\arg\, (\zeta -z)\, |\ <\ \frac{\pi}{4}\ \right\}
\end{equation}
and
\begin{equation}\label{domain}
{\frak {D}}\ =\ \bigcup\limits_{\zeta\in E}\, S_{\zeta}\ \cup D_*
\end{equation}
where
$$
D_*\ =\ \left\{ z\in\mathbb C:\ |z|\ \leq\ \frac{1}{\surd {2}}\
\right\}\ .
$$
It is easy to see that $\partial\frak D$ contains $E$ and is a
rectifiable Jordan curve because $\partial\mathbb D\setminus E$ is
open set and hence the latter consists of a countable collection of
arcs of $\partial\mathbb D$, see the corresponding illustrations in
\cite{Ko}, Section III.D.1.

By the construction, radii of $\mathbb D$ to every $\zeta\in E$
belong to $\frak D$ and the function
$\varphi(\zeta):=\lim\limits_{n\to\infty}\varphi_n(\zeta)$,
$\varphi_n(\zeta):=u(r_n\zeta)$, $n=1,2,\ldots $, with arbitrary
sequence $r_n\to 1-0$ as $n\to\infty$, is measurable, see e.g.
Corollary 2.3.10 in \cite{Fe}. Thus, by the known Egorov theorem,
see e.g. Theorem 2.3.7 in \cite{Fe}, with no loss of generality we
may assume that $\varphi_n\to\varphi$ uniformly on $E$ and that
$\varphi$ is continuous on $E$, see e.g. Section 7.2 in \cite{MRSY}.

Let us consider the sequence of the functions
\begin{equation}\label{functionsD}
\psi_n(\zeta)\ :=\ \sup\limits_{z\in S_{\zeta}\cap\, D^n_{\zeta}}\
|u(z)-\varphi(\zeta)|\ ,\ \ \ \ \ \ \zeta\in E\ ,
\end{equation}
where $D^n_{\zeta}=\{ z\in\mathbb C: |z-\zeta|<\varepsilon_n\}$ with
$\varepsilon_n\searrow 0$ as $n\to\infty$. First of all,
$\psi_n(\zeta)\to 0$ as $n\to\infty$ for every $\zeta\in E$.
Moreover, the functions $\psi_n(\zeta)$ are measurable again by
Corollary 2.3.10 in \cite{Fe} because of
$\psi_n(\zeta)=\lim\limits_{m\to\infty}\, \psi_{mn}(\zeta)$ as
$m\to\infty$ where the functions
\begin{equation}\label{functionsR}
\psi_{mn}(\zeta)\ :=\ \max\limits_{z\in \overline{S_{\zeta}}\cap\,
R^{mn}_{\zeta}}\ |u(z)-\varphi(\zeta)|\ ,\ \ \ \ \ \
R^{mn}_{\zeta}:=\overline{D^n_{\zeta}}\setminus D^{n+m}_{\zeta} \ ,\
\ \ \ \ \ \zeta\in E\ ,
\end{equation}
are continuous. Indeed, $\psi_{mn}(\zeta)$ coincide with the
Hausdorff distance between the compact sets $u(\,
\overline{S_{\zeta}}\cap\, R^{mn}_{\zeta})$ and
$\{\varphi(\zeta)\}$, see e.g. Theorem 2.21.VII in \cite{Ku}, and
any distance is continuous with respect to its variables, recall
also that both functions $u$ and $\varphi$ are continuous.

Again by the Egorov theorem, with no loss of generality we may
consider that $\psi_n\to 0$ uniformly on $E$. The latter implies
that the restriction $U$ of the harmonic function $u$ to the domain
$\frak D$ is bounded. Indeed, let us assume that there exists a
sequence of points $z_n\in\frak D$ such that $|u(z_n)|\ge n$,
$n=1,2,\ldots$. With no loss of generality we may consider that
$z_n\to\zeta\in E$ because the function $u$ is bounded on the
compact subsets of $\mathbb D$ and by the construction $E=\partial
\frak D\cap\partial\mathbb D$ and $E$ is compact. Moreover, by the
construction of $\frak D$, we also may consider that $z_n\in
S_{\zeta_n}$, $\zeta_n\in E$, $n=1,2,\ldots$ and that
$\zeta_n\to\zeta$ as $n\to\infty$. Consequently, it should be that
$u(z_n)\to \varphi(\zeta)$ because $\psi_n(\zeta_n)\to 0$ as
$\zeta_n\to\zeta$, see e.g. Theorem 7.1(2) and Proposition 7.1 in
\cite{MRSY}. The latter conclusion contradicts to the above
assumption.

Further, by the construction the domain $\frak D$ is simply
connected and hence by the Riemann theorem there exists a conformal
mapping $w=\omega(z)$ of $\frak D$ onto $\mathbb D$, see e.g.
Theorem II.2.1 in \cite{Go}. Note that the function $U_*:=U\circ
\omega^{-1}$ is a bounded harmonic function in $\mathbb D$ and there
exists its conjugate harmonic function $V_*$ in $\mathbb D$, i.e.
$F:=U_*+i\, V_*$ is an analytic function in $\mathbb D$. Let $N$ be
a positive number that is greater than $\sup\limits_{w\in\mathbb
D}\, |U_*(w)|= \sup\limits_{z\in\frak D}\, |U(z)|$. Then the
analytic function $g(w):=F(w)/(N-F(w))$, $w\in\mathbb D$, is
bounded. Thus, by the Fatou theorem, see e.g. Corollary III.A in
\cite{Ko}, $g$ has finite angular limits as $w\to W$ for a.e.
$W\in\partial\mathbb D$. By Remark 2 these limits cannot be equal to
1 on a subset of $\partial\mathbb D$ of a positive measure.
Consequently, the function $F(w)$ has also (finite !) angular limits
as $w\to W$ for a.e. $W\in\partial\mathbb D$.

Let us consider the analytic function $f=F\circ\omega$ given in the
domain $\frak D$. By the construction ${\rm Re}\, f=U=u|_{\frak D}$
and hence $V:={\rm Im}\, f$ is its conjugate harmonic function in
$\frak D$. By the standard uniqueness theorem for analytic
functions, we have that $V=v|_{\frak D}$ where $v$ is a conjugate
harmonic function for $u$ in $\mathbb D$. Recall that the latter is
unique up to an additive constant. Thus, it remains to prove that
the function $f(z)$ has (finite !) angular limits as $z\to\zeta$ for
a.e. $\zeta\in E$. For this goal, note that the rectifiable curve
$\partial\frak D$ has tangent a.e. with respect to its natural
parameter. It is clear that tangents at points $\zeta\in E$ to
$\partial\frak D$ (where they exist !) coincide with the
corresponding tangents at $\zeta$ to $\partial\mathbb D$.

By the Caratheodory theorem $\omega$ can be extended to a
homeomorphism of $\overline {\frak D}$ onto $\overline{\mathbb D}$
and, since $\partial {\frak D}$ is rectifiable, by the theorem of F.
and M. Riesz $\mathrm{length}\ \omega^{-1}({\cal E})=0$ whenever
${\cal E}\subset\partial\mathbb D$ with $\mathrm{length}\ {\cal
E}=0$, see e.g. Theorems II.C.1 and II.D.2 in \cite{Ko}. By the
Lindel\"of theorem, see e.g. Theorem II.C.2 in \cite{Ko}, if
$\partial {\frak D}$ has a tangent at a point $\zeta$, then $$\arg\
[\omega(\zeta)-\omega(z)]-\arg\ [\zeta-z]\to{\mathrm const}\ \ \ \ \
{\mathrm as}\ \ \  z\to\zeta\ .$$ In other words, the conformal
images of sectors in $\frak D$ with a vertex at
$\zeta\in\partial\frak D$ is asymptotically the same as sectors in
$\mathbb D$ with a vertex at $w=\omega(\zeta)\in\partial\mathbb D$
up to the corresponding shifts and rotations. Consequently,
nontangential paths in $\mathbb D$ are transformed under
$\omega^{-1}$ into nontangential paths in $\frak D$ and inversely at
the corresponding points of $\partial\mathbb D$ and $\partial\frak
D$.

Thus, in particular, $v(z)$ has finite angular limits
$\varphi_*(\zeta)$ for a.e. $\zeta\in E$. Moreover, the function
$\varphi_*:E\to\mathbb R$ is measurable because
$\varphi_*(\zeta)\,=\, \lim\limits_{n\to\infty}\, v_n(\zeta)$ where
$v_n(\zeta):=v(r_n\zeta)$, $n=1,2,\ldots $, with $r_n\to 1-0$ as
$n\to\infty$, see e.g. Corollary 2.3.10 in \cite{Fe}.  $\Box$

\bigskip

In particular, we have the following consequence of Theorem 1.

\bigskip

{\bf Corollary 1.} {\it Let $u:\mathbb D\to\mathbb R$ be a harmonic
function that has angular limits a.e. on the unit circle
$\partial\mathbb D$. Then its conjugate harmonic functions $v$ in
$\mathbb D$ also have angular limits a.e. on $\partial\mathbb D$ and
both boundary functions are measurable.}

\bigskip

By Remark 1 we have also the next consequence of Theorem 1.

\bigskip

{\bf Corollary 2.} {\it Let $u:\mathbb D\to\mathbb R$ be a harmonic
function in the class $h^1$. Then its conjugate harmonic functions
$v:\mathbb D\to\mathbb R$ have (finite !) angular limits $v(z)\to
\varphi(\zeta)$ as $z\to\zeta$ for a.e. $\zeta\in\partial\mathbb
D.$}

%\bigskip

\cc
\section{The case of rectifiable Jordan domains}

%\bigskip

{\bf Theorem 2.} {\it Let $D$ be a Jordan domain in $\mathbb C$ with
a rectifiable boundary and $u: D\to\mathbb R$ be a harmonic function
that has angular limits on a measurable set $E$ of $\partial D$ with
respect to the natural parameter. Then its conjugate harmonic
functions $v: D\to\mathbb R$ also have (finite !) angular limits
a.e. on $E$ with respect to the natural parameter and both boundary
functions are measurable on $E$ with respect to this parameter.}

\bigskip

{\it Proof.} Again by the Riemann theorem there exists a conformal
mapping $w=\omega(z)$ of $D$ onto $\mathbb D$ and by the
Caratheodory theorem $\omega$ can be extended to a homeomorphism of
$\overline { D}$ onto $\overline{\mathbb D}$. As known, a
rectifiable curves have tangent a.e. with respect to the natural
parameter. Hence $\partial D$ has a tangent at every point $\zeta$
of the set $ E$ except its subset $\cal E$ with $\mathrm{length}\
{\cal E}=0$. By the Lindel\"of theorem, for every $\zeta\in
E\setminus {\cal E}$,
$$\arg\ [\omega(\zeta)-\omega(z)]-\arg\ [\zeta-z]\to{\mathrm const}\
\ \ \ \ {\mathrm as}\ \ \  z\to\zeta\ .$$ Thus, the harmonic
function $u_*:=u\circ \omega^{-1}$ given in $\mathbb D$ has angular
limits $\varphi_*(w)$ at all points $w$ of the set
$E_*:=\omega(E\setminus {\cal E})\subseteq \partial\mathbb D$.
Consequently, by Theorem 1 its conjugate harmonic function
$v_*:\mathbb D\to\mathbb R$ has (finite !) angular limits
$\psi_*(w)$ at a.e. point $w\in E_*$ and the boundary functions
$\varphi_*:E_*\to\mathbb R$ and $\psi_*:E_*\to\mathbb R$ are
measurable. The harmonic function $v:=v_*\circ \omega$ is conjugate
for $u$ because the function $f:=f_*\circ\omega$, where
$f_*:=u_*+v_*$, is analytic. Finally, by theorems of Lindel\"of and
F. and M. Riesz $v$ has (finite !) angular limits $\psi(\zeta)
=\psi_*(\omega(\zeta))$ at a.e. point $\zeta\in E$.

The boundary functions $\varphi =\varphi_*\circ\omega$ and $\psi
=\psi_*\circ\omega$ of $u$ and $v$ on $E$, correspondingly, are
measurable functions on $E$ because
$\varphi(\zeta)=\lim\limits_{n\to\infty}\varphi_n(\zeta)$ for all
$\zeta\in E$ and $\psi(\zeta)=\lim\limits_{n\to\infty}\psi_n(\zeta)$
for a.e. $\zeta\in E$, where the functions
$\varphi_n(\zeta):=u_*(r_n\omega(\zeta))$ and
$\psi_n(\zeta):=v_*(r_n\omega(\zeta))$ with $r_n\to 1-0$ as
$n\to\infty$ are continuous, see e.g. Corollary 2.3.10 in \cite{Fe}.
$\Box$

\bigskip

{\bf Corollary 3.} {\it Let $D$ be a Jordan domain in $\mathbb C$
with a rectifiable boundary and $u: D\to\mathbb R$ be a harmonic
function that has angular limits a.e. on $\partial D$ with respect
to the natural parameter. Then its conjugate harmonic functions $v:
D\to\mathbb R$ also have (finite !) angular limits a.e. on
$\partial D$ and both boundary functions are measurable on $E$ with
respect to the natural parameter.}

\bigskip

{\bf Remark 3.} These results can be extended to domains whose
boundaries consist of a finite number of mutually disjoint
rectifiable Jordan curves  (through splitting into a finite
collection of Jordan's domains !).

\bigskip

The established facts can be applied to various boundary value
problems for harmonic and analytic functions in the plane, see e.g.
\cite{R1}--\cite{R4}.

\medskip
\noindent
{\bf Vladimir Ryazanov,}\\
Institute of Applied Mathematics and Mechanics,\\
National Academy of Sciences of Ukraine, room 417,\\
19 General Batyuk Str., Slavyansk, 84116, Ukraine,\\
vl.ryazanov1@gmail.com

\end{document}